\documentclass[reqno, 11pt]{amsart} 
\usepackage{amssymb,graphicx,bbm}
\usepackage[usenames,dvipsnames]{color}
\usepackage{ifpdf}
\usepackage[pdftex]{thumbpdf}       

\ifpdf
\usepackage[pdftex,colorlinks]{hyperref} 
\hypersetup{%
pdftitle={An alternating service problem},
pdfauthor={M. Vlasiou, I.J.B.F. Adan},
pdfkeywords={},
bookmarksnumbered,
pdfstartview={FitH},
urlcolor=cyan,
}%
\else
  \newcommand\phantomsection\relax
\fi

\newcommand{\m}[1]{\mathcal{#1}}
\newcommand{\e}{\mathbb{E}}
\newcommand{\p}{\mathbb{P}}

\newtheorem{theorem}{Theorem}
\newtheorem{corollary}{Corollary}
\newtheorem{lemma}{Lemma}
\theoremstyle{remark} 
\newtheorem{remark}{Remark}
\newcommand{\lta}[1][]{\mbox{$\alpha^{#1}$}}
\numberwithin{equation}{section}

\newcommand{\Dfb}{\mbox{$F_B$}} 
\newcommand{\Dfa}{\mbox{$F_A$}}
\newcommand{\Dfr}{\mbox{$F_R$}}
\newcommand{\Dfw}{\mbox{$F_W$}} 
\newcommand{\dfw}{\mbox{$f_W$}}

\begin{document}
\title{An alternating service problem}
\author[M. Vlasiou, I.J.B.F. Adan]{M. Vlasiou$^*$, I.J.B.F. Adan$^{*,\,**}$}
\thanks{$^*$ Email address: vlasiou@eurandom.tue.nl}\thanks{$^{**}$Email address: iadan@win.tue.nl}
\date{January 27, 2005}
\maketitle

\begin{center}
$^*$ EURANDOM,\\
P.O. Box 513, 5600 MB Eindhoven, The Netherlands.
\end{center}

\begin{center}
$^{**}$ Eindhoven University of Technology,
\\Department of Mathematics \& Computer Science,
\\P.O. Box 513, 5600 MB Eindhoven, The Netherlands.
\end{center}

\begin{abstract}
We consider a system consisting of a server alternating between two service points. At both service points there is an infinite queue of customers that have to undergo a preparation phase before being served. We are interested in the waiting time of the server. The waiting time of the server satisfies an equation very similar to Lindley's equation for the waiting time in the $GI/G/1$ queue. We will analyse this Lindley-type equation under the assumptions that the preparation phase follows a phase-type distribution while the service times have a general distribution. If we relax the condition that the server alternates between the service points, then the model turns out to be the machine repair problem. Although the latter is a well-known problem, the distribution of the waiting time of the server has not been studied yet. We shall derive this distribution under the same setting and we shall compare the two models numerically. As expected, the waiting time of the server is on average smaller in the machine repair problem than in the alternating service system, but they are not stochastically ordered.
\end{abstract}

\section{Introduction}
In this paper we shall study a model that involves one server alternating between two service points. The model applies in many real-life situations and it is described by a Lindley-type equation. This equation is identical to the original Lindley equation apart from a plus sign that is changed into a minus sign. To better illustrate the model, we give a simple example.

Consider an ophthalmologist who performs laser surgeries for cataracts. Since the procedure lasts only 10 minutes and is rather simple, he will typically schedule many consecutive surgeries in one day. Before surgery, the patient undergoes a preparation phase, which does not require the surgeon's attendance. In order to optimise the doctor's utilisation, the following strategy is followed. There are two operating rooms that work non-stop. While the surgeon works in one of them, the next patient is being prepared in the other one. As soon as the surgeon completes one operation, he moves to the other room and a new patient starts his preparation period in the room that has just been emptied.

Apart from the above example, we may think of a hairdresser that has an assistant to help with the preparation of the customers or of a canteen with one employee and two counters that the employee serves in turns. This model arises naturally also in a two-carousel bi-directional storage system, where a picker serves in turns two carousels; see for example Park \textsl{et al}.\ ~\cite{park03} and Vlasiou \textsl{et al}.\ ~\cite{vlasiou04}.

We want to analyse the waiting time of the server, that is, the surgeon in the above example. We assume that apart from the customer that is being served, there is always at least one more customer in the system. In other words, there is always at least one customer in the preparation phase, which means that the server has to wait only because the next customer may not have completed his preparation phase. Furthermore the server is not allowed to serve two consecutive customers at the same service point and must alternate between the service points.

This condition is crucial. If we remove this condition, then the problem turns out to be the classical machine repair problem. In that setting, there is a number of machines working in parallel (two in our situation) and one repairman. As soon as a machine fails, it joins the repair queue in order to be served. The machine repair problem, also known as the computer terminal model (see for example Bertsekas and Gallager ~\cite{bertsekas-DN}) or as the time sharing system (Kleinrock ~\cite{kleinrock-QS2}, Section 4.11) is a well studied problem in the literature. It is one of the key models to describe problems with a finite input population. A fairly extensive analysis of the machine repair problem can be found in Tak\'{a}cs ~\cite{takacs-ITQ}, Chapter 5. In the following we will compare the two models and discuss their performance.

The issue that is usually investigated in the machine repair problem is the waiting time of a machine until it becomes again operational. In our situation though we are concerned with the waiting time of the repairman. This question has not been treated in the classical literature, perhaps because in the machine repair problem the operating time of the machine is usually more valuable than the utilisation of the repairman. In Section \ref{s:comparison} we shall compare the waiting times of the repairman in the classical machine repair problem and our model. We shall show that the random variables for the waiting time in the two situations are not stochastically ordered. However, on average, the alternating strategy leads to longer waiting times for the server. Furthermore we will show that the probability that the server does not have to wait is larger in the alternating service system than in the non-alternating one. This result is perhaps counterintuitive, since the inequality for the mean waiting times of the server in the two situations is reversed.

In the following section we shall introduce the model and explain the interesting aspects of it and the implications in analysis of the different sign in the Lindley-type equation for the waiting time. In Section \ref{s:B_erlang} we shall derive the distribution of the waiting time of the server, provided that the preparation time of a customer follows an Erlang or a phase-type distribution. Continuing with Section \ref{s:MachineRepair}, we shall introduce the machine repair model and analyse the waiting time of the repairman, or in other words we shall remove the restriction that the server alternates between the service points. We will compare the two models in the next section and we will conclude with Section \ref{s:NumericalResults} with some numerical results.

\section{The model}\label{s:model}
We consider a system consisting of one server and two service points. At each service point there is an infinite queue of customers that needs to be served. The server alternates between the service points, serving one customer at a time. Before being served by the server, a customer must undergo first a preparation phase. Thus the server, after having finished serving a customer at one service point, may have to wait for the preparation phase of the customer at the other service point to be completed. We are interested in the waiting time of the server.
Let $B_n$ denote the preparation time for the $n$-th customer and let $A_n$ be the time the server spends on this customer. Then the waiting times $W_n$ of the server satisfy the following relation:
\begin{equation*}
  W_{n+1}=\max\{0, B_{n+1}-A_n-W_n\}.
\end{equation*}
We assume that $\{A_n\}$ and $\{B_n\}$ are independent and identically distributed (i.i.d.)\ sequences of nonnegative random variables, which are mutually independent and have finite means. Further, for every $n$, $A_n$ follows some general distribution $\Dfa(\cdot)$ and $B_n$ follows a phase-type distribution denoted by $\Dfb(\cdot)$. Clearly, the stochastic process $\{W_n\}$ is an aperiodic regenerative process with a finite mean cycle length with the time points where $W_n=0$ being the regeneration points. Therefore, there exists a unique limiting distribution. In the following we shall suppress all subscripts when we refer to the system in equilibrium. Let $W$ be the random variable with this limiting distribution, then
\begin{equation}\label{recursion}
  W\stackrel{\m{D}}{=}\max\{0, B-A-W\},
\end{equation}
where $A$  and $B$ are independent and distributed according to $\Dfa(\cdot)$ and $\Dfb(\cdot)$ respectively.

Note the striking similarity to Lindley's equation for the waiting times in a single server queue. If only the sign of $W_n$ were different, then we would be analysing the waiting time in the $GI/PH/1$ model. Lindley's equation is one of the most studied equations in queueing theory. For excellent textbook treatments we refer to Asmussen ~\cite{asmussen-APQ}, Cohen ~\cite{cohen-SSQ} and the references therein. It is a challenging problem to investigate the implications of this subtle difference between the two equations.

For this model we shall try to obtain an explicit expression for the distribution $\Dfw(\cdot)$ of the waiting time $W$\@. Let $\omega(\cdot)$ denote the Laplace transform of $W$, i.e.,
$$
\omega(s)=\int_0^\infty e^{-sx} d\, \Dfw(x).
$$
The derivative of order $i$ of the transform is $\omega^{(i)}(\cdot)$ and by definition $\omega^{(0)}(\cdot)=\omega(\cdot)$. Similarly we define the Laplace transform $\lta(\cdot)$ of the random variable $A$\@. To keep expressions simple, we shall also use the function $\phi$ defined as $\phi(s)=\omega(s) \, \lta(s)$. We can now proceed with the analysis.

\section{The waiting time distribution}\label{s:B_erlang}
In the following we shall derive the distribution of the waiting time of the server, assuming that the service time $A$ follows some general distribution and the preparation time $B$ follows a phase-type distribution. The phase-type distributions that we will consider are mixtures of Erlang distributions with the same scale parameters. Therefore we will first consider the case where $B$ follows an Erlang-$n$ distribution, that we denote by $E_n(\cdot)$.

\subsection{Erlang preparation times}
Let $B$ be the sum of $n$ independent random variables $X_1,\ldots,X_n$ that are exponentially distributed with parameter $\mu$. The use of Laplace transforms is a standard approach for the analysis of Lindley's equation. Hence it is natural to try this approach for equation \eqref{recursion}. Then we can readily prove the following.

\begin{theorem}[Alternating service system] \label{th:density}
The waiting time distribution has a mass $p_0$ at the origin, which is given by
$$
p_0=\p[B<W+A]=1-\sum_{i=0}^{n-1}\frac{(-\mu)^i}{i!}\,\phi^{(i)}(\mu)
$$
and has a density $\dfw(\cdot)$ on $[0, \infty)$ that is given by
\begin{equation}\label{A density}
   \dfw(x)=\mu^ne^{-\mu x}\sum_{i=0}^{n-1}\frac{(-1)^i}{i!}\,\phi^{(i)}(\mu)\frac{x^{n-1-i}}{(n-1-i)!}.
\end{equation}
In the above expression
$$
\phi^{(i)}(\mu)=\sum_{k=0}^i \binom{i}{k} \omega^{(k)}(\mu)\,\lta[(i-k)](\mu)
$$
and the parameters $\omega^{(i)}(\mu)$ for $i=0,\ldots,n-1$ are the unique solution to the system of equations
\begin{align}\label{balance eq}
\omega(\mu)&=1-\sum_{i=0}^{n-1}(-\mu)^i(1-\frac{1}{2^{n-i}}) \sum_{k=0}^i \frac{\omega^{(k)}(\mu)\,\lta[(i-k)](\mu)}{k!\,(i-k)!}\\ \notag
\omega^{(\ell)}(\mu)&=\sum_{i=0}^{n-1}\mu^{i-\ell}\frac{(-1)^{i+\ell}}{2^{n-i+\ell}}\frac{(n-i+\ell-1)!}{(n-i-1)!}\sum_{k=0}^i \frac{\omega^{(k)}(\mu)\,\lta[(i-k)](\mu)}{k!\,(i-k)!}\quad\mbox{for $\ell=1,\ldots,n-1$.}
\end{align}
\end{theorem}

\begin{proof}
We use the following notation: $\e[X ; A]=\e[X \cdot \mathbbm{1}_{[A]}]$. Consider the Laplace transform of \eqref{recursion}; then we have that
\begin{align}
\nonumber   \omega(s)=\e[e^{-sW}]&=\p[B<W+A]+\e[e^{-s(B-W-A)} \,; B\geqslant W+A]\\
        \label{eq:sum} &=\p[B<W+A]+\e[e^{-s(B-W-A)} \,; X_1\geqslant W+A]\\
        \nonumber      &\qquad+\sum_{i=1}^{n-1}\e[e^{-s(B-W-A)} ; X_1+\cdots+X_i\leqslant W+A\leqslant X_1+\cdots+X_{i+1}]
\end{align}
Using standard techniques and the memoryless property of the exponential distribution, one can show that
\begin{align}
 \nonumber \e[e^{-s(B-W-A)} ; &X_1\geqslant W+A]\\
          \nonumber           &=\e[e^{-s(X_2+X_3+\ldots+X_n)}e^{-s(X_1-W-A)} ; X_1\geqslant W+A]\\
                   \nonumber  &=\left(\frac{\mu}{\mu+s}\right)^{n-1}\e[e^{-s(X_1-W-A)} ; X_1\geqslant W+A]\\
                    \nonumber &=\left(\frac{\mu}{\mu+s}\right)^{n-1}\e[e^{-s(X_1-W-A)} \mid X_1\geqslant W+A]\,\p[X_1\geqslant W+A]\\
                     \nonumber&=\left(\frac{\mu}{\mu+s}\right)^{n}\p[X_1\geqslant W+A]\quad\mbox{(due to the memoryless property)}\\
                     \label{a}&=\left(\frac{\mu}{\mu+s}\right)^{n}\omega(\mu)\lta(\mu).
\end{align}
Additionally, for $Y_i=X_1+\cdots+X_i$ we have that
\begin{multline}
\label{b}  \e[e^{-s(B-W-A)} ; Y_i\leqslant W+A\leqslant Y_{i+1}]=\\
  \left(\frac{\mu}{\mu+s}\right)^{n-i-1}\e[e^{-s(Y_{i+1}-W-A)} \mid Y_i\leqslant W+A\leqslant Y_{i+1}]\,\p[Y_i\leqslant W+A\leqslant Y_{i+1}]=\\
 \left(\frac{\mu}{\mu+s}\right)^{n-i}\frac{(-\mu)^i\phi^{(i)}(\mu)}{i!}.
\end{multline}

Finally, we calculate the probability $\p[B<W+A]$ by substituting $s=0$ in \eqref{eq:sum} and using equations \eqref{a} and \eqref{b}. Straightforward calculations give us now that
\begin{equation}\label{eq:transform_Erl}
  \omega(s)=1-\sum_{i=0}^{n-1}\frac{(-\mu)^i}{i!}\phi^{(i)}(\mu)\left(1-\left(\frac{\mu}{\mu+s}\right)^{n-i}\right).
\end{equation}
Inverting the transform yields the density \eqref{A density}.

Furthermore, the terms $\omega^{(i)}(\mu)$, $i=0,\ldots,n-1$, that are included in $\phi^{(i)}(\mu)$ still need to be determined. To obtain the values of  $\omega^{(i)}(\mu)$, for $i=0,\ldots,n-1$, we differentiate \eqref{eq:transform_Erl} $n-1$ times and we evaluate $\omega^{(i)}(s)$,  $i=0,\ldots,n-1$ at the point $s=\mu$. This gives us the system of equations \eqref{balance eq}. The fact that the solution of the system is unique follows from the general theory of Markov chains that implies that there is a unique equilibrium distribution and thus also a unique solution to \eqref{balance eq}.
\end{proof}

\begin{corollary}
The throughput $\theta$ satisfies
$$
\theta^{-1}= \e[W]+ \e[A]=\sum_{i=0}^{n-1}\frac{(-1)^i}{i!}\,\phi^{(i)}(\mu) \mu^{i-1}(n-i)-\lta['](0).
$$
\end{corollary}

It is quite interesting to note that the density of the waiting time can be rewritten as
$$
\dfw(x)=\mu e^{-\mu x}\sum_{i=1}^n p_i \frac{(\mu x)^{i-1}}{(i-1)!},
$$
where
$$
p_i=\frac{(-\mu)^{n-i}\phi^{(n-i)}(\mu)}{(n-i)!}
$$
is the probability that directly after a service completion exactly $i$ exponential phases of $B$ remain.

\begin{figure}[hbp]
\begin{center}
\includegraphics[width=0.8\textwidth]{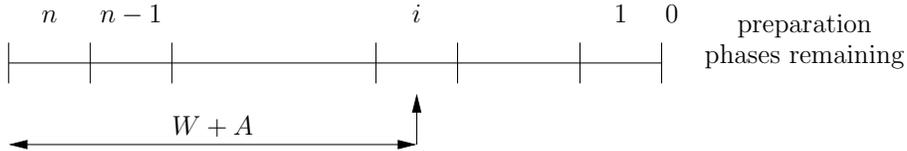}
\end{center}
\caption{The waiting time has a mixed Erlang distribution.}
\label{fig:explanation}
\end{figure}

As it is also clear from Figure \ref{fig:explanation}, with probability $p_i$ the distribution of the waiting time is Erlang-$i$, for $i=1, \ldots, n$. Furthermore, the probability $p_0$ that the server does not have to wait, or equivalently that at least $n$ exponential phases expired, is $p_0=1-\sum_{i=1}^n p_i$.

So practically the problem is reduced to obtaining the solution of an $n \times n$ linear system. Extending the above result to mixtures of Erlang distributions is simple.

\subsection{Phase-type preparation times}
For $n=1,\ldots,N$, let the random variable $X_n$ follow an Erlang-$n$ distribution with parameter $\mu$ and let the random variable $B$ of the preparation times be equal to $X_n$ with a probability $\kappa_n$. In other words the distribution function of $B$ is given by
\begin{equation}\label{mixed}
\Dfb(x)=\sum_{n=1}^N \kappa_n\left(1-e^{-\mu x}\sum_{j=0}^{n-1}\frac{(\mu x)^j}{j!}\right) ,\qquad x \geqslant 0.
\end{equation}
This class of phase-type distributions may be used to approximate any given distribution for the preparation times arbitrarily close; see Schassberger ~\cite{schassberger-W}. Below we show that Theorem \ref{th:density} can be extended to service distributions of the form \eqref{mixed}.

By conditioning on the number of phases of $B$, we find that
$$
\omega(s)=\e[e^{-sW}]=\p[B<W+A]+\sum_{n=1}^N \kappa_n \e[e^{-s(X_n-W-A)} ; X_n\geqslant W+A].
$$
Since $X_n$ follows now an Erlang-$n$ distribution, the last equation is practically a linear combination of equation \eqref{eq:sum}, summed over all probabilities $\kappa_n$ for $n=1,\ldots,N$. This means that we can directly use the analysis of Section \ref{s:B_erlang} to calculate the Laplace transform of $W$ in this situation (cf.\ equation \eqref{eq:transform_Erl}). So we have that
\begin{equation}\label{whatever}
\omega(s)=1-\sum_{n=1}^N \kappa_n \sum_{i=0}^{n-1}\frac{(-\mu)^i}{i!}\phi^{(i)}(\mu)\left(1-\left(\frac{\mu}{\mu+s}\right)^{n-i}\right),
\end{equation}
where the terms $\phi^{(i)}(\mu)$ can be calculated in a similar fashion as previously. Inverting \eqref{whatever} yields the density of the following theorem (cf.\ Theorem \ref{th:density}).

\begin{theorem}
Let \eqref{mixed} be the distribution of the random variable $B$\@. Then the distribution of the server's waiting time has mass $p_0$ at zero which is given by
$$
p_0=\p[B<W+A]=1-\sum_{n=1}^N \sum_{i=0}^{n-1}\kappa_n\frac{(-\mu)^i}{i!}\,\phi^{(i)}(\mu)
$$
and has a density on $[0, \infty)$ that is given by
\begin{equation*}
  \dfw(x)=\sum_{n=1}^N \kappa_n\left(\mu^ne^{-\mu x}\sum_{i=0}^{n-1}\frac{(-1)^i}{i!}\,\phi^{(i)}(\mu)\frac{x^{n-1-i}}{(n-1-i)!}\right).
\end{equation*}
\end{theorem}

One can already see from the above theorem the effect of the different sign in the Lindley-type equation that describes our model. The waiting time distribution for the $GI/PH/1$ queue is a mixture of exponentials with different scale parameters (cf.\ Adan and Zhao ~\cite{adan96}). In our case we have that the waiting time distribution is a mixture of Erlang distributions with the same scale parameter for all exponential phases.

As we have mentioned before, the practice of alternating between the service points is inevitably followed in many situations. Still it seems reasonable to argue that it would be more efficient to choose to serve the first customer that has completed his preparation time. If we drop the assumption that the server alternates between the service points then we have the classical machine repair problem.

\section{The machine repair problem}\label{s:MachineRepair}
In the machine repair problem there is a number of machines that are served by a unique repairman when they fail. The machines are working independently and as soon as a machine fails, it joins a queue formed in front of the repairman where it is served in order of arrival. A machine that is repaired is assumed to be as good as new. The model described in Section \ref{s:model} is exactly the machine repair problem after we drop the assumption that the server alternates between the service points. There are two machines that work in parallel (the two service points), the preparation time of the customer is equivalent to the life time of the machine until it fails and the service time of the customer is the time the repairman needs to repair the machine.

What we are interested in is the waiting time of the repairman until a machine breaks down or, in other words, the waiting time of the server until the preparation phase of one of the customers is completed. It is quite surprising that although the machine repair problem under general assumptions is thoroughly treated in the literature, this question remains unanswered. We would like to compare the two models and to this end we first need to derive the distribution of the waiting time of the server, when the system is in steady state. In the following we will refer to the server or customers instead of the repairman or machines in order to illustrate the analogies between the two models.

Let $B$ be the random variable of the time needed for the preparation phase and $R$ be the remaining preparation time just after a service has been completed. Then obviously the waiting time of the server is $W=\min\{B, R\}$. The random variables $B$ and $R$ are independent, so in order to calculate the distribution of $W$ we need the distribution of $R$\@. In agreement to the alternating service model, $B$ follows an Erlang-$n$ distribution. Note that we do not have a simple Lindley-type recursion for $W$ and therefore this system cannot be easily treated with Laplace transforms. That means that we have to try an alternative approach.

The system can be fully described by the number of remaining phases of preparation time that a customer has to complete, immediately after a service completion. The state space is finite, since there can be at most $n$ phases remaining and the Markov chain is aperiodic and irreducible, so there is a unique equilibrium distribution $\{\pi_i, i= 0,\ldots,n\}$. After completing a service, the other customer may be already waiting for the server (so the $n$ exponential phases of the Erlang-$n$ distribution of the preparation have expired) or he is during one of the $n$ phases of the preparation time. That means that the remaining preparation time $R$ that the server sees immediately after completing a service follows the mixed Erlang distribution $\Dfr(x)=\pi_0+\pi_1 E_1(x)+\ldots+\pi_n E_n(x)$.

So in order to derive the distribution of $R$ (and consequently the distribution of $W$), we need to solve the equilibrium equations $\pi_i=\sum_k \pi_k\, p_{ki}$, $i=0,\ldots,n$, in conjunction with the normalising equation $\sum_k \pi_k=1$, where $p_{ki}$ are the one-step transition probabilities. Let us determine the probabilities $p_{ij}$, for all $i, j\in\{0,\ldots,n\}$.

A transition from state $i$ to state $j$, for $i, j\in\{1,\ldots,n\}$, can be achieved in two ways: either the customer that has just been served or the other one will finish the preparation phase first. Suppose that the customer that has just been served finishes first. In that case we know that the last event just before the service starts is that the $n$-th phase of that customer expired. The other customer was in state $k$ and  during the service time the other customer reached state $j$, i.e.\ $k-j$ phases of that customer have expired. The probability of this event is given by
$$
\sum_{k=j}^n \left(\frac{1}{2}\right)^{n+i-k}\,\binom{n+i-k-1}{n-1} \frac{(-\mu)^{k-j}}{(k-j)!}\,\lta[(k-j)](\mu),
$$
where $\lta(\cdot)$ is as before the Laplace transform of the service time. Note that in the above expression we have that
$$
\p[\mbox{exactly $k-j$ Exp$(\mu)$ phases expired during $[0,A)$}]=\frac{(-\mu)^{k-j}}{(k-j)!}\,\lta[(k-j)](\mu).
$$

Similarly we can determine the probability of a transition from state $i$ to state $j$ in the second case. So in the end we have that for all $i, j\in\{1,\ldots,n\}$,
\begin{equation}\label{1}
p_{ij}=\sum_{k=j}^n \left(\frac{1}{2}\right)^{n+i-k}\,\left[\binom{n+i-k-1}{n-1}+\binom{n+i-k-1}{n-k}\right] \frac{(-\mu)^{k-j}}{(k-j)!}\,\lta[(k-j)](\mu).
\end{equation}

The transition probabilities from state zero to any state to state $i=1,\ldots,n$ are
\begin{equation}\label{2}
p_{0i}=\frac{(-\mu)^{n-i}}{(n-i)!}\,\lta[(n-i)](\mu),
\end{equation}
since starting from state zero means that the other customer was already waiting when the repairman finished a service and reaching state $i$ means that during the service time, exactly $n-i$ exponential phases expired. For the transition from state zero to state zero we have that during the service time at least $n$ exponential phases expired, so
\begin{equation}\label{3}
p_{00}=\sum_{i=n}^\infty \frac{(-\mu)^{i}}{(i)!}\,\lta[(i)](\mu).
\end{equation}
Similarly, we have that for $i=1,\ldots,n$
\begin{equation}\label{4}
p_{i0}=\sum_{k=1}^n \left(\frac{1}{2}\right)^{n+i-k}\,\left[\binom{n+i-k-1}{n-1}+\binom{n+i-k-1}{n-k}\right] \left(\sum_{j=k}^\infty \frac{(-\mu)^{j}}{(j)!}\,\lta[(j)](\mu)\right),
\end{equation}
where $\binom{a}{b}=0$ for $0\leqslant a<b$.

With the one-step transition probabilities one can determine the equilibrium distribution and thus $\Dfr(\cdot)$. Then we have that the distribution of the waiting time of the server, if we drop the assumption that he is alternating between the service points, is given by the following theorem.
\begin{theorem}[Non-alternating service system]
The waiting time distribution is
$$
\Dfw(x)= \Dfr(x) +\Dfb(x)-\Dfr(x)\Dfb(x),
$$
where $\Dfr(\cdot)$ is the distribution of the remaining preparation time of a customer and is equal to
$$
\Dfr(x)=\pi_0+\pi_1 E_1(x)+\ldots+\pi_n E_n(x).
$$
In the above expression $\{\pi_i, i= 0,\ldots,n\}$ is the unique solution to the system of equations
$$
\pi_i=\sum_{k=0}^{n} \pi_k\, p_{ki} \quad\mbox{and } \sum_{k=0}^{n} \pi_k=1, \quad \mbox{for }i=0,\ldots,n,
$$
where $p_{ij}$ are given by the equations \eqref{1}-\eqref{4}.
\end{theorem}

\begin{remark}
The above results can be easily extended to phase-type preparation times of the form \eqref{mixed}. However this extension does not contribute significantly to the analysis, since it is along the same lines of the analysis in this section.
\end{remark}
This method of defining a Markov chain through the remaining phases of the preparation time after a service has been completed and using the equilibrium distribution in order to calculate the mixing probabilities of $R$ can, of course, also be used in the alternating service system. In that case, the waiting time $W$ is exactly the remaining preparation time $R$. Then the probabilities $p_i$, for $i=0,\ldots,n$ as defined in Section \ref{s:B_erlang} will be the equilibrium distribution of the underlying Markov chain. Furthermore, the system of equations \eqref{balance eq} can be rewritten as follows:
\begin{align*}
\omega(\mu)&=p_0+\sum_{i=1}^{n}\frac{p_i}{2^i}\\
\omega^{(\ell)}(\mu)&=\sum_{i=1}^{n}\frac{p_i (-\mu)^{-\ell}}{2^{i+\ell}}\frac{(i+\ell-1)!}{(i-1)!}\quad\mbox{for $\ell=1,\ldots,n-1$.}
\end{align*}

In the next session we shall study various performance characteristics of the two systems.

\section{ Performance comparison}\label{s:comparison}
One may wonder if there is any connection between the waiting time of the server in the two models that can help in understanding how the models perform. From this point on we will use the superscript $\mathrm{A}$ ($\mathrm{NA}$) for all variables associated with the (non-)alternating service system when we specifically want to distinguish between the two situations. Otherwise the superscript will be suppressed. So, for example, the random variable $W^{\mathrm{A}}$ will be the waiting time of the server in the alternating service system.

\subsection{Stochastic ordering}\label{ss:st order}
Suppose that the distributions of the two random variables $X$ and $Y$ have a common support. Then the stochastic ordering $X \geqslant_{st} Y$ is defined as (cf.\ \cite{lehmann55,mann-whitney47,szekli-SODAP})
$$
\p [X\geqslant x]\geqslant \p [Y\geqslant x], \quad\mbox{for all $x$ in the support,}
$$
and we say that $X$ {\it dominates} $Y$.

\begin{figure}[htbp]
\leavevmode
\begin{center}
\includegraphics[width=0.6\textwidth]{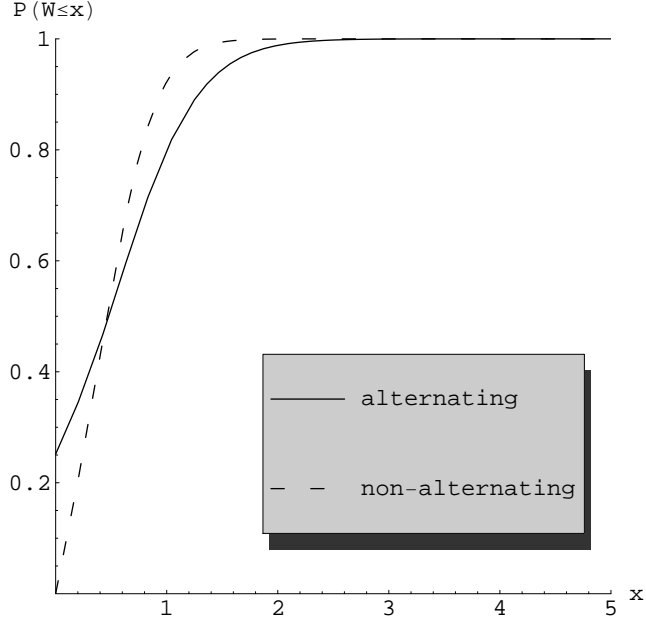}
\end{center}
\caption{$W^{\mathrm{A}}$ and $W^{\mathrm{NA}}$ are not stochastically ordered.}
\label{fig:distr}
\end{figure}

Intuitively one may argue that $W^{\mathrm{A}}\geqslant_{st} W^{\mathrm{NA}}$ since one expects that large waiting times occur with higher probability in the alternating service system. However this is not true. Let us imagine the situation where the service times are equal to zero. Then in the alternating service system we will have that the waiting time of the server is zero if $B_i \geqslant B_{i+1}$ for some $i$. So since $\p[B_i \geqslant B_{i+1}]>0$, we have $\p [W^{\mathrm{A}}=0]>0$. In the non-alternating system however, we will have zero waiting time only if both preparation phases finish at exactly the same instant. Since the preparation times are continuous random variables we have that $\p[B_i=B_{i+1}]=0$ for every $i$ and thus $\p [W^{\mathrm{NA}}=0]=0$. In Figure \ref{fig:distr} we have plotted the distribution of the waiting time for both situations in the case where the service times are equal to zero and the $B$ follows an Erlang-$5$ distribution with $\mu=5$.

The situation that we have described above is not a rare example. In fact, the following result holds.

\begin{theorem}\label{th:W=0}
For any distribution of the preparation and the service time, we have that $\p[W^{\mathrm{A}}=0]\geqslant\p[W^{\mathrm{NA}}=0]$.
\end{theorem}

\begin{proof}
Both processes regenerate when a zero waiting time of the server occurs. Therefore in a cycle there is precisely one customer for whom the server did not have to wait. This means that the fraction of customers for whom the server does not wait is
$$
\p[W=0]=\frac{1}{\e[N]},
$$
where $\e[N]$ is the average number of customers in a cycle, i.e.\ the mean cycle length. So it suffices to show that $\e[N^{\mathrm{A}}]\leqslant \e[N^{\mathrm{NA}}]$.

To prove this, we will couple the two systems and use sample path arguments. We will show that for a given initial state and for any realisation of preparation and service times the number of customers in a cycle is greater in the alternating case than in the non-alternating case. To couple the systems we will use the same realisations for the preparation and the service times. To this end, let $\{{B}_i\}$ be a sequence of preparation times and $\{{A}_i\}$ a sequence of service times. We need to observe the system until the completion of the first cycle. For both systems assume that the server starts servicing the first customer at time zero while at the other service point a customer has just started his preparation phase ${B}_1$. Additionally , let ${R}_n$ be the remaining preparation time for the $n$-th customer, immediately after the service of the $n-1$ customer has finished. As long as ${R}_n\leqslant {B}_{n+1}$ both processes are identical, since both servers will alternate between the two service points. In addition, all waiting times until that point will be strictly positive. As soon as ${R}_n>{B}_{n+1}$, the alternating service system will regenerate for the first time, since we will have that ${W}_n^{\mathrm{A}}={R}_n$ and ${W}_{n+1}^{\mathrm{A}}=0$. The non-alternating system however does not necessarily regenerate. For this system we have that ${W}_n^{\mathrm{NA}}={B}_{n+1}$ and ${R}_{n+1}^{\mathrm{NA}}={R}_n-{B}_{n+1}-{A}_n$. Therefore, if ${R}_{n+1}^{\mathrm{NA}}=0$ then ${W}_{n+1}^{\mathrm{NA}}=0$ and both processes regenerate. Otherwise the non-alternating system will not regenerate. Hence, for each realisation we have that $N^{\mathrm{NA}}\geqslant N^{\mathrm{A}}$ which implies that the mean cycle length $\e[N^{\mathrm{NA}}]$ of the non-alternating system is at least as long as the mean cycle length $\e[N^{\mathrm{A}}]$ of the alternating system.
\end{proof}

\subsection{Mean waiting times}
Although the waiting times in the two situations are not stochastically ordered, we have however that the mean waiting time of the server of the alternating service model $\e[W^{\mathrm{A}}]$ is larger than or equal to the mean waiting time of the server in the non-alternating system $\e[W^{\mathrm{NA}}]$. This is quite natural, since we expect the non-alternating system to perform better in terms of throughput, regardless of the distribution of the preparation phase.

To prove this result for the mean waiting times, we will again couple the two systems. We will make use of the same realisations $\{B_i\}$ and $\{A_i\}$ for the preparation and the service times respectively and we will continue with sample path arguments. We assume that the initial conditions for both systems are the same, i.e.\ at time zero the server starts servicing the first customer, while at the other service point a customer has just started his preparation phase. Then, for the alternating service system, define:\\
$\mbox{$\quad$} D^{\mathrm{A}}_i$: the $i-$th departure time\\
$\mbox{$\quad$} H^{\mathrm{A}}_i$: the time the server can start serving the other service point after time $D^{\mathrm{A}}_i$.\\
Also define in the same way $D^{\mathrm{NA}}_i$ and $H^{\mathrm{NA}}_i$ for the non-alternating system.  We need the following lemma.

\begin{lemma}\label{lemma}
For  all $i$, we have that $D^{\mathrm{A}}_i \geqslant D^{\mathrm{NA}}_i$ and $H^{\mathrm{A}}_i \geqslant H^{\mathrm{NA}}_i$.
\end{lemma}

\begin{proof}
We will apply induction. For $i=1$ we have that $D^{\mathrm{A}}_1\geqslant D^{\mathrm{NA}}_1$ and $H^{\mathrm{A}}_1\geqslant H^{\mathrm{NA}}_1$, since
$$
D^{\mathrm{A}}_1=A_1\geqslant A_1=D^{\mathrm{NA}}_1
$$
and thus
$$
H^{\mathrm{A}}_1=\max\{D^{\mathrm{A}}_1,B_1\} \geqslant \max\{D^{\mathrm{NA}}_1, B_1\}=H^{\mathrm{NA}}_1.
$$

Suppose that for some $i$ we have that $D^{\mathrm{A}}_{i-1} \geqslant D^{\mathrm{NA}}_{i-1}$ and $H^{\mathrm{A}}_{i-1} \geqslant H^{\mathrm{NA}}_{i-1}$. We will prove that $D^{\mathrm{A}}_i \geqslant D^{\mathrm{NA}}_i$ and $H^{\mathrm{A}}_i \geqslant H^{\mathrm{NA}}_i$ and this will conclude the proof.

The first relation is obvious. From the induction hypothesis we have $H^{\mathrm{A}}_{i-1} \geqslant H^{\mathrm{NA}}_{i-1}$, so
$$
D^{\mathrm{A}}_i=H^{\mathrm{A}}_{i-1}+A_i \geqslant \min\{H^{\mathrm{NA}}_{i-1} , D^{\mathrm{NA}}_{i-1}+B_{i}\}+A_i=D^{\mathrm{NA}}_i.
$$

For the second inequality, first notice that
$$
H^{\mathrm{A}}_{i}=\max\{D^{\mathrm{A}}_{i}, D^{\mathrm{A}}_{i-1}+B_{i}\}
\mbox{ and }
H^{\mathrm{NA}}_{i}=\max\{D^{\mathrm{NA}}_{i} , \max\{H^{\mathrm{NA}}_{i-1} , D^{\mathrm{NA}}_{i-1}+B_{i}\}\},
$$
because, for example, in the non-alternating case the other service point will either be ready at time $D^{\mathrm{NA}}_{i}$ when the previous customer departs, or it will be ready after the preparation phase at this point is completed, at the time point equal to the maximum of $H^{\mathrm{NA}}_{i-1}$ and $D^{\mathrm{NA}}_{i-1}+B_{i}$.

To prove that
\begin{equation}\label{h2 inequality}
H^{\mathrm{A}}_{i}=\max\{D^{\mathrm{A}}_{i}, D^{\mathrm{A}}_{i-1}+B_{i}\} \geqslant \max\{D^{\mathrm{NA}}_{i} , \max\{H^{\mathrm{NA}}_{i-1} , D^{\mathrm{NA}}_{i-1}+B_{i}\}\}=H^{\mathrm{NA}}_{i},
\end{equation}
we will show that the maximum term of the left hand side of the inequality \eqref{h2 inequality} is greater than or equal to any term of the right hand side, thus also greater than or equal to the maximum of them.

Assume that $H^{\mathrm{A}}_i=D^{\mathrm{A}}_i$. Then
$D^{\mathrm{A}}_i \geqslant D^{\mathrm{NA}}_i$ as we have proven above, furthermore
$D^{\mathrm{A}}_i=H^{\mathrm{A}}_{i-1}+A_i \geqslant H^{\mathrm{NA}}_{i-1}$
        since $H^{\mathrm{A}}_{i-1} \geqslant H^{\mathrm{NA}}_{i-1}$
        and finally since $H^{\mathrm{A}}_i=D^{\mathrm{A}}_i$ then
$D^{\mathrm{A}}_i \geqslant D^{\mathrm{A}}_{i-1}+B_{i} \geqslant D^{\mathrm{NA}}_{i-1}+B_{i}$.
The case for $H^{\mathrm{A}}_i=D^{\mathrm{A}}_{i-1}+B_{i}$ follows similarly.
\end{proof}

A corollary of the previous result is the following.
\begin{corollary}
For all $i$, $\sum_j^i W^{\mathrm{A}}_j \geqslant_{st} \sum_j^i W^{\mathrm{NA}}_j$.
\end{corollary}
\begin{proof}
  The proof is a direct consequence of the fact that for the coupled systems
$$
W^{\mathrm{A}}_1+A_1+\ldots+W^{\mathrm{A}}_i+A_i = D^{\mathrm{A}}_i \geqslant D^{\mathrm{NA}}_i = W^{\mathrm{NA}}_1+A_1+\ldots+W^{\mathrm{NA}}_i+A_i.
$$
\end{proof}
So, although the random variables $W^{\mathrm{A}}$ and $W^{\mathrm{NA}}$ are not stochastically ordered, the partial sums of the sequences ${W^{\mathrm{A}}_i}$ and ${W^{\mathrm{NA}}_i}$ are.

It is also interesting to note that Lemma \ref{lemma} immediately implies that the throughput is greater in the non-alternating system than in the alternating system since
$$
\theta^{\mathrm{A}}=\lim_{i\rightarrow \infty} \frac{i}{D^{\mathrm{A}}_i}\leqslant\lim_{i\rightarrow \infty} \frac{i}{D^{\mathrm{NA}}_i}=\theta^{\mathrm{NA}}.
$$
Moreover we have that $\theta= \left(\e[W]+ \e[A]\right)^{-1}$, so we can readily establish the following result:

\begin{theorem}\label{th:mean waiting}
Given any distribution for the preparation and the service time, we have that  $\e[W^\mathrm{A}] \geqslant \e[W^{\mathrm{NA}}]$.
\end{theorem}

\begin{figure}[ht]
\leavevmode
\begin{center}
\includegraphics[width=0.8\textwidth]{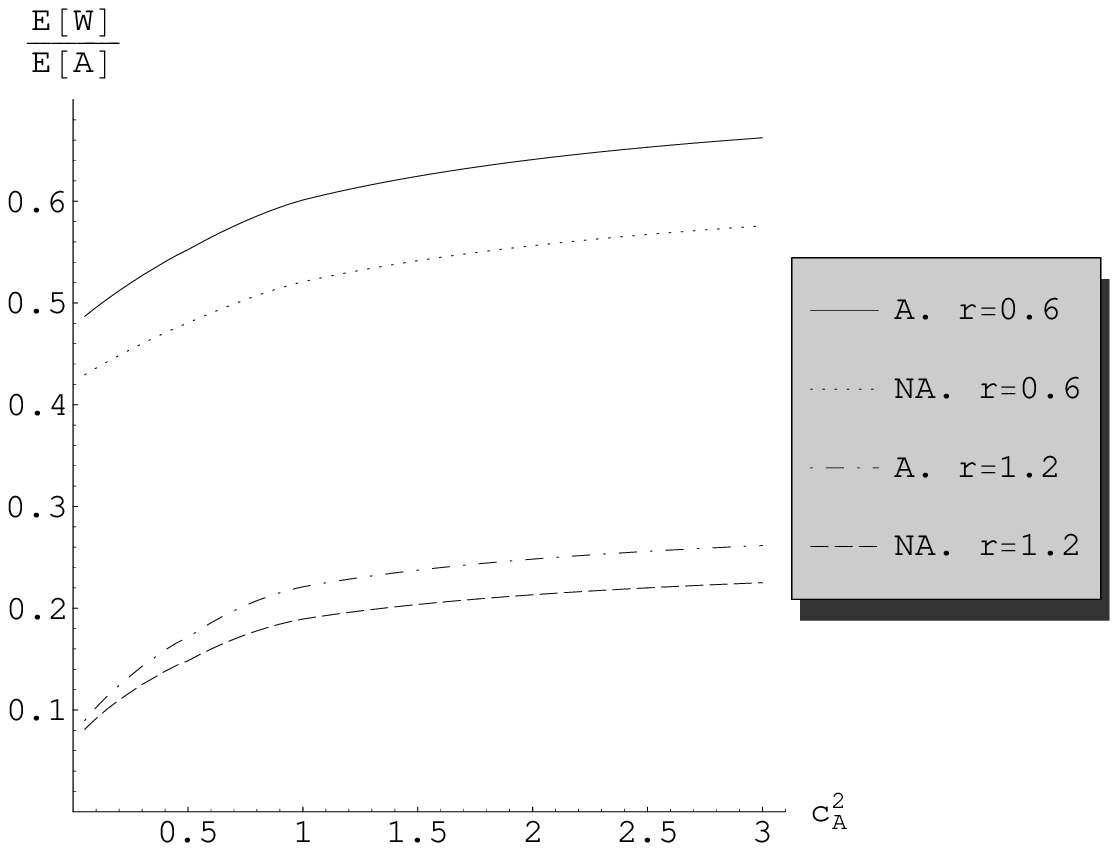}
\end{center}
\caption{$\e[W^\textrm{A}]$ is greater than or equal to $\e[W^\textrm{NA}]$.}
\label{fig:EW-cv2pick}
\end{figure}

Figure \ref{fig:EW-cv2pick} demonstrates a typical situation. For two values of the ratio $r=\e[A]/\e[B]$ we have plotted the normalised waiting time $\e[W]/\e[A]$ versus the squared coefficient of variation $c^2_A$ of the service time $A$\@. We have chosen the mean service time to be $\e[A]=1$ and the preparation time to be composed of five exponential phases. As before, \textrm{A} stands for the alternating service system and \textrm{NA} for the non-alternating system. One can see from these two examples that the average waiting time in the alternating service system is larger than in the non-alternating system. As it is the case for the $GI/G/1$ queue, the waiting time depends almost linearly on $c^2_A$.  As $c^2_A$ increases, the waiting time also increases and for the alternating case the rate of change is bigger. The difference of the mean waiting time in the alternating and the non-alternating case is eventually almost constant and this difference increases as the value of $r$ decreases. In Appendix \ref{s:recipe} we give more details on the way we chose the distribution for the service time.

\begin{remark}
From Theorem \ref{th:W=0} and Theorem \ref{th:mean waiting} we can conclude that there is at least one point where the waiting time distributions of both systems intersect. Figure \ref{fig:distr} suggests though that this point is unique. So, since the mean waiting times are both finite, this implies that $W^{\mathrm{NA}}$ is smaller than $W^{\mathrm{A}}$ with respect to the \emph{increasing convex ordering}; namely
$$
\e [\phi(W^{\mathrm{NA}})]\leqslant\e [\phi(W^{\mathrm{A}})]
$$
for all increasing convex functions $\phi$, for which the mean exists.
This follows as a direct application of the Karlin-Novikoff cut-criterion (cf.\ Szekli ~\cite{szekli-SODAP}).
\end{remark}

\section{Numerical results}\label{s:NumericalResults}
This section is devoted to some numerical results. In Figure \ref{fig:EW-cv2pick} we have already shown how the normalised waiting time changes when the squared coefficient of variation of the service time is modified.

\begin{figure}[htbp]
\includegraphics[width=0.45\textwidth]{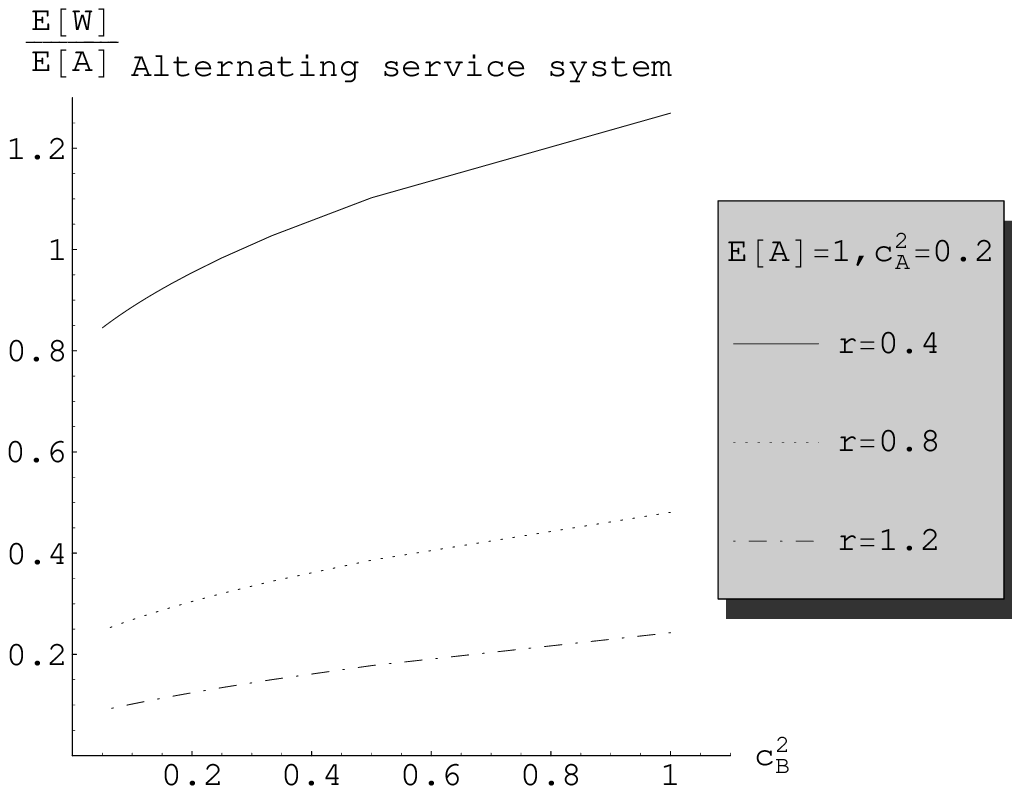}
\hfill
\includegraphics[width=0.45\textwidth]{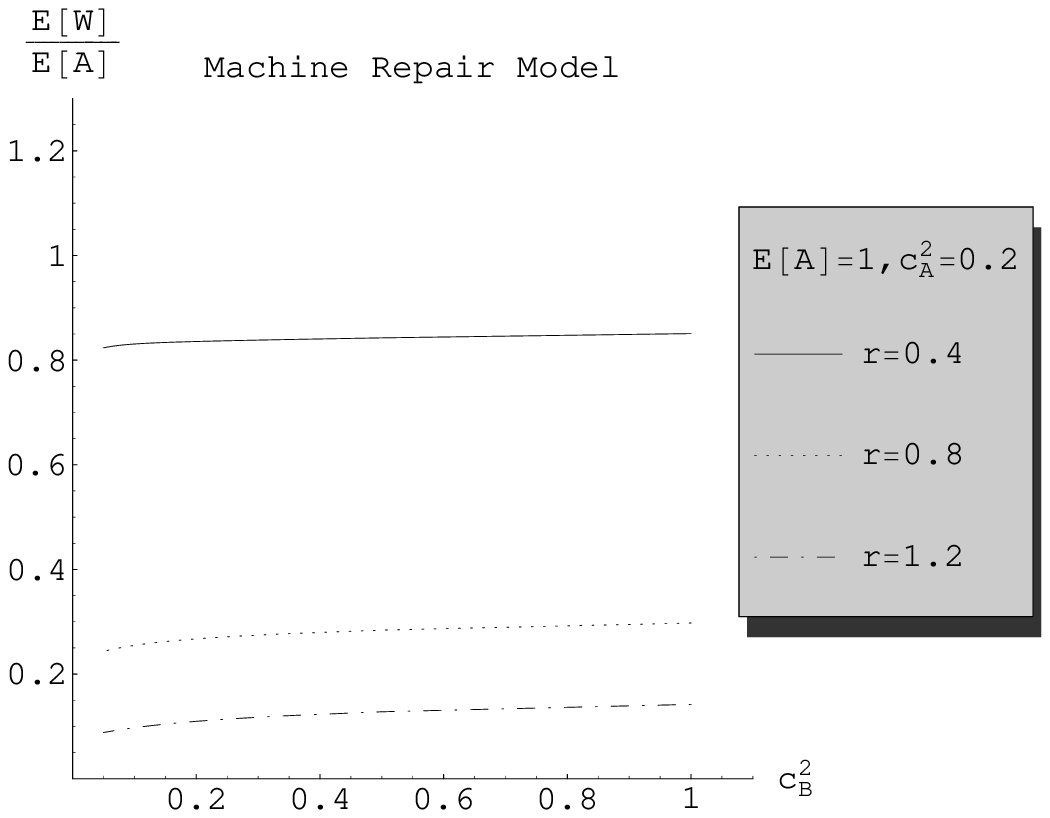}
\caption{The normalised waiting time is almost insensitive to $c^2_B$ in the non-alternating system.}
\label{fig:EW-cv2rot}
\end{figure}

Figure \ref{fig:EW-cv2rot} shows the normalised waiting time plotted against the squared coefficient of variation of the preparation time. The preparation time is assumed to follow an Erlang distribution. We chose  $\e[A]=1$ and $c_A^2=0.2$ and we fitted a mixed Erlang distribution according to the procedure described in Appendix \ref{s:recipe}\@. We have plotted the normalised waiting time for three different values of the ratio $r$; namely for $r=0.4$, which implies that the service time is 40\% of the preparation time, up to $r=1.2$. The latter implies that for the alternating service model the server in general does not have to wait much. One can see that the normalised waiting time depends almost linearly on $c^2_B$ for the alternating service system, but for the non-alternating system it is almost insensitive to $c^2_B$ and thus to the number of exponential phases of the preparation time. This can be explained by the fact that Erlang loss models are insensitive to the service time distribution apart from its first moment; see for example Kelly ~\cite{kelly-RSN}. More specifically, one can view the machine repair model that we have described, as an $E_n/G/2/2$ loss system. Here the repairman would act as the Poisson source of an Erlang loss model if $B$ would follow an exponential distribution. However the preparation times are a sum of exponentials and that causes the slight fluctuation in the mean waiting time.

\begin{figure}[htbp]
\includegraphics[width=0.45\textwidth]{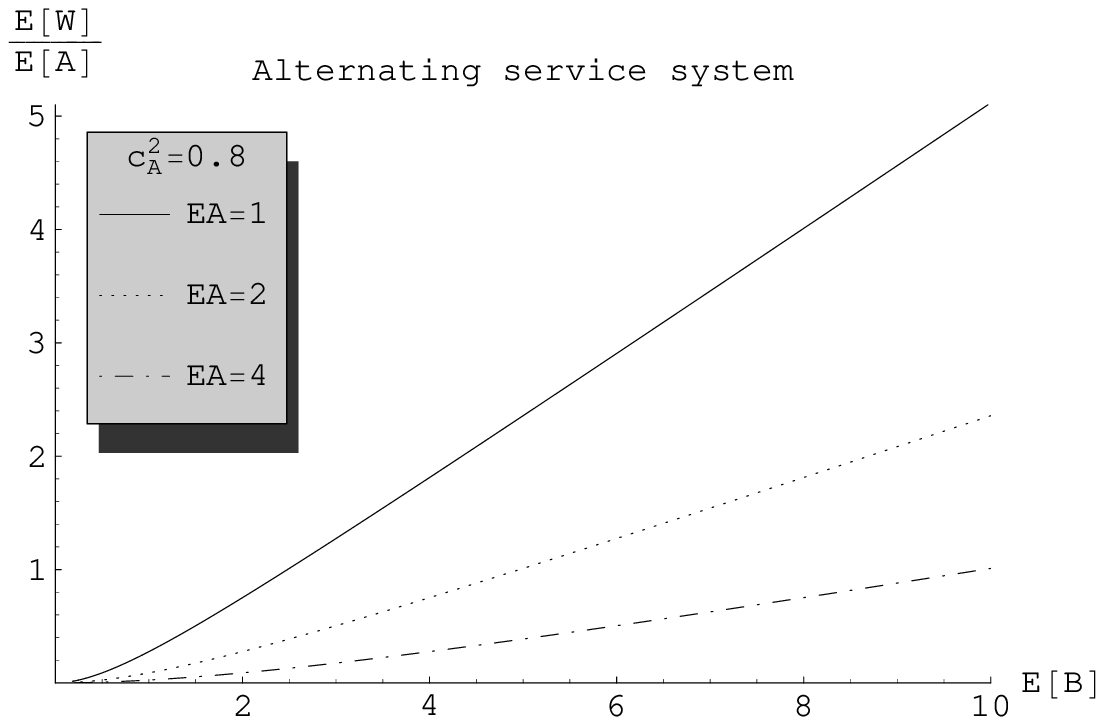}
\hfill
\includegraphics[width=0.45\textwidth]{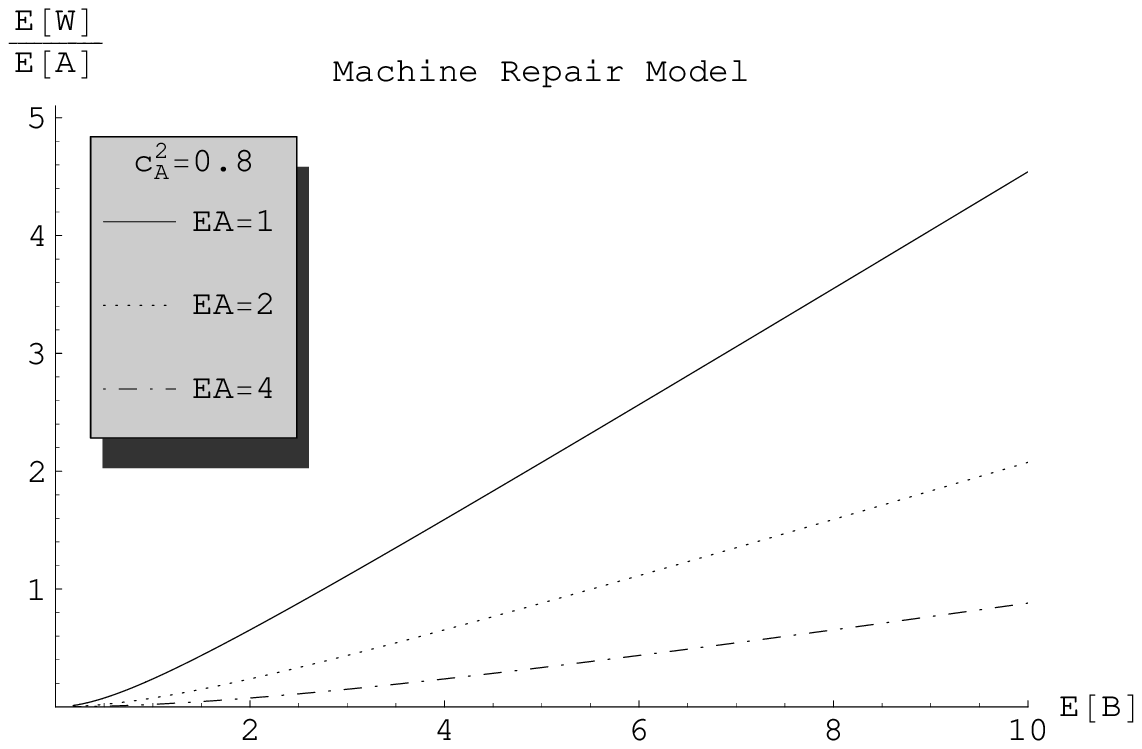}
\caption{The normalised waiting time vs. $\e[B]$.}
\label{fig:EW-Erot}
\end{figure}

Figure \ref{fig:EW-Erot} shows the normalised waiting time plotted against the mean preparation time. We have chosen $c^2_A$ to be equal to 0.8 and we have fitted a mixed Erlang distribution to the mean service time and the squared coefficient of service. As expected the normalised waiting time $\e[W]/\e[A]$ depends almost linearly on the mean preparation time. For larger values of the mean preparation time, the normalised waiting time increases.

\appendix
\section{Fitting distributions}\label{s:recipe}
In Figure \ref{fig:EW-cv2pick} we chose the mean service time $\e[A]$ to be equal to one and we plotted the normalised waiting time versus $c^2_A$. For each setting we fitted a mixed Erlang or hyperexponential distribution to $\e[A]$ and $c^2_A$, depending on whether the squared coefficient of variation is less or greater than $1$ (see, for example, Tijms ~\cite{tijms-SM}). More specifically, if $1/n \leqslant c_A^2 \leqslant 1/(n-1)$ for some $n = 2, 3, \ldots$, then the mean and squared coefficient of variation of the mixed Erlang distribution
$$
\Dfa(x)= p
\left(1-e^{-\mu x}\sum_{j=0}^{n-2}\frac{(\mu x)^j}{j!}\right)
+ (1-p)
\left(1-e^{-\mu x}\sum_{j=0}^{n-1}\frac{(\mu x)^j}{j!}\right),
\qquad x \geqslant 0 ,
$$
matches with $\e[A]$ and $c_A^2$, provided the parameters $p$ and $\mu$ are chosen as
$$
p = \frac{1}{1+c_A^2}
[ n c_A^2 - \{ n(1+c_A^2) - n^2 c_A^2 \}^{1/2} ] , \qquad \mu = \frac{n-p}{\e[A]} .
$$
On the other hand, if $c_A^2 > 1$, then the mean and squared coefficient of variation of the hyperexponential distribution
$$
\Dfa(x) = p_1 (1- e^{-\mu_1 x}) + p_2 (1-  e^{-\mu_2 x}) ,
\qquad x \geqslant 0,
$$
match with $\e[A]$ and $c_A^2$ provided the parameters $\mu_1, \mu_2, p_1$ and $p_2$ are chosen as
\begin{eqnarray*}
&& p_1 =
\frac{1}{2} \left( 1 + \sqrt{\frac{c_A^2 - 1}{c_A^2 + 1}}\right) , \qquad p_2 = 1-p_1, \\
&& \mu_1 =\frac{2 p_1}{\e[A]} \qquad \mbox{and}\qquad\mu_2 = \frac{2 p_2}{\e[A]}.
\end{eqnarray*}

\bibliographystyle{apt}
\bibliography{maria}

\end{document}